\documentclass[a4paper,12pt]{article}

\usepackage{amssymb}

\usepackage{amssymb}
\usepackage{amsmath}
\usepackage{amsfonts}                                                           
\usepackage{mathrsfs} 

\usepackage{bbm}

\usepackage{graphicx} 
\usepackage{float} 

\usepackage{arydshln}

\usepackage[fleqn,tbtags]{mathtools}

\def \C{\hbox{\sf \rlap{\kern.25em \vrule width.1em height1.6ex depth-.1ex}C}}
\def \D{{\sf I\kern-1.5ptD\,}}
\def \K{{\sf I\kern-1.5ptK\,}}
\def \N{{\sf I\kern-1.5ptN\,}}
\def \P{{\sf I\kern-1.5ptP\,}}
\def \Q{\hbox{\sf \rlap{\kern.25em \vrule width.1em height1.6ex depth-.1ex}Q}}
\def \R{{\sf I\kern-1.5ptR\,}}
\def \T{{\sf T\kern-6.5ptT\,}}
\def \Z{{\sf Z\kern-5.0ptZ\,}}

\def\begitem#1 {\bigskip\pagebreak[1]%
     \refstepcounter{subsection}{\nopagebreak[4]%
     \thesubsection\hskip 0.5truecm}
     {\sc#1}\hskip 1pt.\nopagebreak[4]\par\nopagebreak[4]%
      \begin{enumerate}\rm\nopagebreak[4]}
\def\BEGITEMK#1 #2{\bigskip\pagebreak[1]%
      \refstepcounter{subsection}{\nopagebreak[4]%
     \thesubsection\hskip 0.5truecm}\nopagebreak[4]
     {\bf#1}\hskip 1pt.\nopagebreak[4]\par\nopagebreak[4]%
     \medskip\nopagebreak[4]\rm#2\nopagebreak[4]%
     \begin{enumerate}\nopagebreak[4]\rm}
\def\enditem{\end{enumerate}}

\newcommand{\qed}{\hfill \ensuremath{\Box}}

\newtheorem{theorem}{Theorem}[section]
\newtheorem{corollary}{Corollary}[section]
\newtheorem{lemma}{Lemma}[section]
\newtheorem{proposition}{Proposition}[section]

\usepackage{centernot}

\newcommand{\atom}{\operatorname*{atom}}
\newcommand{\Cay}{\operatorname{Cay}}
\newcommand{\ICG}{\operatorname{ICG}}

\begin{document}

\title{Adding generators in cyclic groups}
\author{{\bf J.W.~Sander\footnote{Corresponding author}}\\
Institut f\"ur Mathematik und Angewandte Informatik,
Universit\"at Hildesheim,\\
D-31141 Hildesheim, Germany\\
{\tt sander@imai.uni-hildesheim.de}\\[.1in]
and\\[.1in]
{\bf T.~Sander}\\ 
Fakult\"at f\"ur Informatik,
Ostfalia Hochschule f\"ur angewandte Wissenschaften,\\
D-38302 Wolfenb\"uttel, Germany\\
{\tt t.sander@ostfalia.de}}
\maketitle

\begin{abstract}

For a cyclic group $\langle a\rangle$, define the atom of $a$ as the set of all elements generating $\langle a\rangle$. 
Given any two elements $a,b$ of a finite cyclic group $G$,
we study the sumset of the atom of $a$ and the atom of $b$.
It is known that such a sumset is a disjoint union of atoms. The goal of this paper is 
to offer a deeper understanding of this phenomenon, by 
determining which atoms make up the sum of two given atoms and by computing
the exact number of representations of each element of the sumset.

\end{abstract}


\section{Introduction}

When adding multiplicative objects, as for instance primes (cf.~Goldbach's conjecture) or units in a ring, one usually is rather sceptical to find a lot of algebraic structure. Let $\mathbb Z_n:= \mathbb Z/n\mathbb Z$ denote the ring of residue classes $\bmod\, n$ and $\mathbb Z_n^*$ its group of units, i.e.~the multiplicative group of primitive residues $a \bmod n$ with $(a,n)=1$. Now consider the sumset $\mathbb Z_n^*+\mathbb Z_n^*$, where generally  $A+B:=\{a+b:\; a\in A,\, b\in B\}$ for any two non-empty sets $A$ and $B$ such that the addition makes sense. 

In the following, let $(x,y)$ denote the greatest common divisor of $x$ and $y$, let $\varphi$
denote Euler's totient function, and let $\mathbb P_{\geq 3}$ be the set of odd primes.
In 2000 {\sc Deaconescu} \cite{DEA} derived from earlier work together with {\sc Du} \cite{DEA2} on the number of similar automorphisms in finite cyclic groups a formula for 
\begin{equation}   \label{jnt}
N_n(c):=\# \{(x,y)\in (\mathbb Z_n^*)^2: \; x+y\equiv c \bmod n\},
\end{equation}
for arbitrary $c\in \mathbb Z_n$.  It follows from this formula that $N_n(c)>0$ for all $n$ and $c$ except for the case where $n$ is even and $c$ is odd, which obviously has to be excluded. Hence  $\mathbb Z_n^*+\mathbb Z_n^* = \mathbb Z_n$ for all odd $n$. This may not seem too surprising, but rather as a probabilistic than an algebraic phenomenon due to the fact that $\# \mathbb Z_n^* = \varphi(n)$  is large compared with $\# \mathbb Z_n =n$. The same view is taken by most mathematicians with regard to Goldbach's conjecture, reading $\mathbb P_{\geq 3} + \mathbb P_{\geq 3} = \{ 2n: \, n\geq 3\}$ in our context. 

However, in 2009 the first author \cite{SAN} of the present article 
gave a proof for an extension of \eqref{jnt} by using multiplicativity of $N_n(c)$ with respect to $n$. The reader might share the astonishment called forth by observing such a strong structural feature when adding multiplicative objects like units in a ring.
Of course, the addition of residue classes has been studied long before. 
For example, in the 1930ies {\sc Davenport} \cite{DAV} and {\sc Chowla} \cite{CHO}
gave bounds on the size of the sumset of any two sets of residue classes with applications to Waring's problem. In what follows,
we restrict ourselves to special sets of residue classes. The benefit is that we gain a very
precise understanding of the resulting sumset structure. The setting has been inspired 
by a problem in graph theory (see Section \ref{sec:gr}), but it turns out to be an
interesting object of study on its own.

Given a finite abelian group $G$ and an element $a\in G$, define the \textit{atom} of $a$ as
\[   
{\rm atom}(a) := \{a'\in G: \; \langle a' \rangle = \langle a \rangle \}. 
\]
In other words, ${\rm atom}(a)$ is the set of all generators of $\langle a \rangle$. 
Clearly, the atoms of any two elements of $G$ are either disjoint or identical. 
The term ``atom'' originates from the theory of Boolean algebras where it denotes 
the second minimal elements of a lattice.
In this case it refers to the Boolean algebra generated by the subgroups of $G$. It is not difficult to see 
that every element of this Boolean algebra is a disjoint union of atoms. 

Recently {\sc Klotz} and the second author proved the following result in \cite{KLO}: \newline
{\it Let $G$ be a finite abelian group with $a,b\in G$. Then 
${\rm atom}(a) + {\rm atom}(b)$ 
is the (disjoint) union of atoms of $G$. } 


In this article we extend the result of {\sc Klotz} and {\sc T. Sander}  \cite{KLO} in the case where $G$ is a cyclic group. It turns out that in this situation the above phenomenon can be explained much more explicitly.  In particular we shall see the union of which atoms make up the sum of two given atoms (cf.~Theorem \ref{thm2}). Along the way (cf.~Theorem \ref{thm1}) we shall determine the number of representations of each element in ${\rm atom}(a) + {\rm atom}(b)$, which substantially generalizes the formula for the counting function $N_n(c)$ in \eqref{jnt} obtained by {\sc J.W. Sander} \cite{SAN}.


\section{Terminology and the Reduction Lemma}

Up to isomorphism there is exactly one cyclic group of order $n$ for each positive integer $n$, and a standard model is the residue class group $\mathbb Z_n$ with respect to addition $\bmod\, n$. 
We use this model, because in our context 
it is profitable  to consider $\mathbb Z_n$ as a ring, also equipped with multiplication $\bmod\, n$, which turns $\mathbb Z_n$ into a principal ideal domain. 

If we denote by $R^*$ the set of all units in a ring $R$ with $1$, then $R^*$ is a multiplicative group and, in particular,  $\mathbb Z_n^* = \{u\in \mathbb Z_n: \; (u,n)=1\}$. From this point of view, the cyclic subgroup $\langle a \rangle$  of the additive group $\mathbb Z_n$ simply is the principal ideal $(a):=a\mathbb Z_n$ generated by $a$ in the ring $\mathbb Z_n$, and ${\rm atom}(a)$ is exactly the set of all generators of $(a)$, i.e. 
\begin{equation}  \label{atomunit}
{\rm atom}(a)= a\mathbb Z_n^* = \{ au: \; 1\leq u \leq n, \; (u, n)=1\}.
\end{equation}
Usually the set on the righthand side of \eqref{atomunit} provides multiple representations of the elements of ${\rm atom}(a)$. For that reason, we shall prefer the alternative 
\begin{equation}  \label{atom}
{\rm atom}(a)= (a)^{\text{\tiny $\bigstar$}} := 
\{ ax: \; 1\leq x \leq {\rm ord}(a), \; (x, {\rm ord}(a))=1\},
\end{equation}
where ${\rm ord}(a)= \frac{n}{(a,n)}$ denotes the order of $a$ in the additive group $\mathbb Z_n$, and then by \eqref{atom} each element of ${\rm atom}(a)$ is uniquely represented.
A proof for the identity $a\mathbb Z_n^* = (a)^{\text{\tiny $\bigstar$}}$ as well as the uniqueness of the representation in \eqref{atom} is given in Corollary 
\ref{corunit} below. 
We like to draw the reader's attention to the two different types of asterisks we use.
While $R^*$ denotes the unit group in a ring $R$ with $1$, we have $1\notin (a)$ for $(a,n)>1$, and in this case $(a)^{\text{\tiny $\bigstar$}}$ is not the set of units in the ideal $(a)$.
However, all $x\in (a)^{\text{\tiny $\bigstar$}}$ are of type $(a,n)\cdot u$ for some unit $u\in \mathbb Z_n^*$. Observe that our notation includes the definition of $I^{\text{\tiny $\bigstar$}}$ for the zero ideal, namely $\{0\}^{\text{\tiny $\bigstar$}}=\{0\}$.

The theorem of {\sc Klotz} and {\sc T. Sander} \cite{KLO} now reads:  \newline
{\it Let $I$ and $J$ be two ideals in the residue class ring $\mathbb Z_n$. Then $I^{\text{\tiny $\bigstar$}}+J^{\text{\tiny $\bigstar$}}$ is the (disjoint) union of atoms in $\mathbb Z_n$, i.e.~there are ideals $I_1,\ldots,I_k \subseteq \mathbb Z_n$, say, such that
\begin{equation}   \label{sumofatoms} 
I^{\text{\tiny $\bigstar$}}+ J^{\text{\tiny $\bigstar$}} = \bigcup_{j=1}^k \; I_j^{\text{\tiny $\bigstar$}}.  
\end{equation}
}
It is our main goal to determine $I_1,\ldots,I_k$ explicitly in terms of $I$ and $J$.

In order to simplify matters 
we take advantage of the natural order in the set of (positive) integers. Any ideal $I$ in the ring $\mathbb Z_n$ is principal, hence there is some $a\in\mathbb Z_n$ satisfying $(a)=I$. The generating element $a$ is uniquely determined if we require $a$ to be represented by the least non-negative residue $\bmod\, n$ among all generators of $I$. This minimal generator $a$ of $I$ will be called the \textit{leader} ${\rm lead}(I)$ of $I$. Clearly, ${\rm lead}(I)\mid n$ and ${\rm ord}({\rm lead}(I)) = \frac{n}{{\rm lead}(I)}$ for any ideal $I \subseteq \mathbb Z_n$. Moreover, we have by (\ref{atom})
\begin{equation}  \label{idealunits}
I^{\text{\tiny $\bigstar$}}  =  \Big\{ {\rm lead}(I)\cdot x: \; 1\leq x \leq \tfrac{n}{{\rm lead}(I)}\,, \; \big(x,\tfrac{n}{{\rm lead}(I)}\big)=1\Big\}, 
\end{equation}
hence $|I^{\text{\tiny $\bigstar$}}| = \varphi(\frac{n}{{\rm lead}(I)})$.

Given two ideals $I$ and $J$ in $\mathbb Z_n$, we first identify those $c\in \mathbb Z_n$ lying in $I^{\text{\tiny $\bigstar$}}+J^{\text{\tiny $\bigstar$}}$ and determine for each such $c$ the number of representations. Finally, we deduce the desired decomposition \eqref{sumofatoms}.
Therefore, let $a:={\rm lead}(I)$ and $b:={\rm lead}(J)$, hence $a\mid n$ and $b\mid n$,  and define for any $c\in \mathbb Z_n$
\[
S(c)=S_{n;a,b}(c):= \{(u,v)\in I^{\text{\tiny $\bigstar$}}\times J^{\text{\tiny $\bigstar$}}: \; u+v =c\}.
\]
By (\ref{idealunits}), we obtain
\begin{equation}  \label{basicset}
S_{n;a,b}(c)= \{(ax,by): \;1\leq x \leq \tfrac{n}{a}, \,1\leq y \leq \tfrac{n}{b},\,  
(x,\tfrac{n}{a})=(y,\tfrac{n}{b})=1 ,\, ax+by\equiv c\bmod n
\},
\end{equation}
and
\begin{equation}  \label{basic}
N(c)= N_{n;a,b}(c) := \# S_{n;a,b}(c) = \mathop{\sum_{\genfrac{}{}{0pt}{1}{1\leq x \leq \frac{n}{a}}{(x,\frac{n}{a})=1}} 
\sum_{\genfrac{}{}{0pt}{1}{1\leq y \leq \frac{n}{b}}{(y,\frac{n}{b})=1}}}\limits_{ax+by\equiv c\bmod n} 1.
\end{equation}

It will facilitate further considerations if we may assume the leaders 
of the ideals $I$ and $J$ to be coprime. This is justified by the following lemma proved in Section \ref{technical}.
\begin{lemma} [Reduction Lemma]  \label{redlemma}
Let $n$ be a positive integer with divisors $a$ and $b$ and $g:=(a,b)$. For any $c\in \mathbb Z_n$ satisfying $g\mid c$, and on setting
$n':=\frac{n}{g}$, $a':=\frac{a}{g}$, $b':= \frac{b}{g}$ and $c':=\frac{c}{g}$, we have:
\begin{itemize}
\item[(i)] The function
\[
\rho_c: \left\{ \begin{array}{ccc}
             S_{n;a,b}(c) & \rightarrow & S_{n';a',b'}(c') \\
             (ax,by)    &  \mapsto        & (a'x,b'y)
             \end{array}  \right.
\]
is 1-1.
\item[(ii)]
$N_{n;a,b}(c)= N_{n';a',b'}(c')$. 
\end{itemize}
\end{lemma}

{\sc Proof of Lemma \ref{redlemma}. }
Since (ii) is an immediate consequence of (i), it suffices to show (i). Since $g\mid c$, the numbers $n',a',b',c'$ are integers. Let $(ax,by)\in S_{n;a,b}(c)$, hence $1\leq x \leq \frac{n}{a}$ and 
$1\leq y \leq \frac{n}{b}$ with $(x,\frac{n}{a})=(y,\frac{n}{b})=1$ by \eqref{basicset}. Moreover, we have $ax+by\equiv c \bmod n$, which is equivalent with $a'x+b'y\equiv c' \bmod n'$. By the fact that $\frac{n}{a}=\frac{n'}{a'}$ and $\frac{n}{b}=\frac{n'}{b'}$, we conclude that $(a'x,b'y)\in S_{n';a',b'}(c')$, i.e.~$\rho_c$ is well defined. It is also obvious that $\rho_c$ is $1-1$.

\qed \newline
The Reduction Lemma tells us in case $(a,b)>1$ how to obtain all representations of $c$ in $S_{n;a,b}(c)$ from the representations of $c'$ in $S_{n';a',b'}(c')$.

In order to be able to evaluate the double sum \eqref{basic} in a satisfactory manner, we introduce some more terminology.
For positive integers $m$ and $k$ we define
\[
\varphi^*(m,k) := m \prod_{\genfrac{}{}{0pt}{1}{p\in \mathbb P}{p\mid m,\, p\mid k}} \left(1-\frac{1}{p}\right) \prod_{\genfrac{}{}{0pt}{1}{p\in \mathbb P}{p\mid m,\, p\nmid k}} \left(1-\frac{2}{p}\right). 
\]
This modified version of Euler's totient function, which for fixed $k$ is multiplicative with respect to $m$, was introduced by the first author in \cite{SAN} (in slightly different notation). For any positive integer $m$, we denote by ${\rm rad}(m):= \prod_{p\in \mathbb P,\, p\mid m} p\;$ the so-called \textit{radical} or \textit{squarefree kernel} of $m$.


\section{Main results}

Theorem \ref{thm1} below is a generalisation of Theorem 1.1 in \cite{SAN}, where the first author proved that
\[
N_{n;1,1}(c)=
\# \{ (x,y)\in \mathbb Z_n^*\times \mathbb Z_n^*: \; x+y \equiv c \bmod n\} = \varphi^*(n,c).
\]
This is the special case $a=b=1$ of the following result.
\begin{theorem}  \label{thm1}
Let $n$ be a positive integer with divisors $a$ and $b$ and $g:=(a,b)$, and let $c\in \mathbb Z_n$.
\begin{itemize}
\item[(i)] If $g\nmid c$, then $N_{n;a,b}(c) = 0$.
\item[(ii)] Let $g\mid c$ and set $n':=\frac{n}{g}$, $a':=\frac{a}{g}$, $b':= \frac{b}{g}$ and $c':= \frac{c}{g}$. If $(c',a'b')>1$, then $N_{n;a,b}(c) = 0$. If $(c',a'b')=1$, then
\begin{align}  \label{genform}
\begin{split}
N_{n;a,b}(c) 
&= m \prod_{\genfrac{}{}{0pt}{1}{p\in \mathbb P}{p\mid m, \,p\mid a'b'}}
\left(1-\frac{1}{p}\right)
\prod_{\genfrac{}{}{0pt}{1}{p\in \mathbb P}{p\mid n',\, p\nmid a'b', \, p\mid c'}} \left(1-\frac{1}{p}\right)
\prod_{\genfrac{}{}{0pt}{1}{p\in \mathbb P}{p\mid n',\,  p\nmid a'b'c'}} \left(1-\frac{2}{p}\right)\\
&= \frac{m}{{\rm rad}(m)}\, \varphi(m_1)\,\varphi(m_2)\,\varphi^*(m_3,c'),  
\end{split}
\end{align}
where $m:=\frac{n'}{a'b'}$ and we write ${\rm rad}(m) = m_1 m_2m_3$ with $m_1\mid a'$, $m_2\mid b'$ and $(m_3,a'b')=1$. 
\end{itemize}
\end{theorem}

\begin{corollary}  \label{cor1}
Let $n$ be a positive integer with divisors $a$ and $b$, $g:=(a,b)$, and let $c\in \mathbb Z_n$. Then $c\in (a)^{\text{\tiny $\bigstar$}} + (b)^{\text{\tiny $\bigstar$}}$ 
 if and only if the following three conditions are satisfied:
\[
(i)\;\; g\mid c;\quad (ii)\;\; (c',a'b')=1; \quad (iii)\;\; n'
\mbox{ is odd or } a'b'c' \mbox{ is even},
\]
where $n':=\frac{n}{g}$, $a':=\frac{a}{g}$, $b':= \frac{b}{g}$ and $c':= \frac{c}{g}$.
\end{corollary}

\begin{corollary}  \label{cor2}
Let $n$ be a positive integer with divisors $a$ and $b$, and let $I$ be an ideal in $\mathbb Z_n$. Then $N_{n;a,b}(u)=N_{n;a,b}(v)$ for any $u,v\in I^{\text{\tiny $\bigstar$}}$. 
\end{corollary}

\begin{theorem}  \label{thm2}
Let $n$ be a positive integer with divisors $a$ and $b$ and $g:=(a,b)$,
and let $c\in \mathbb Z_n$. We set $n':=\frac{n}{g}$, $a':=\frac{a}{g}$ and $b':= \frac{b}{g}$.
\begin{itemize}
\item[(A)] If $2\nmid n'$ or $2\mid a'b'$, we have
\begin{itemize}
\item[(1)] $N_{n;a,b}(c)>0$ if and only if $g\mid c$ and $(c',a'b')=1$ for $c':= \frac{c}{g}$.
\item[(2)] If $c\neq 0$ and $N_{n;a,b}(c)>0$, then there is a unique ideal $I\subset \mathbb Z_{n'}$ such that $c'\in I^{\text{\tiny $\bigstar$}}$. Moreover, $\mbox{lead}(I)= (c',\widetilde{m}_3)$, where $\widetilde{m}_3$ is the largest divisor of $\frac{n'}{a'b'}$ satisfying $(\widetilde{m}_3,a'b')=1\;$. 
\item[(3)] {$\displaystyle (a)^{\text{\tiny $\bigstar$}} + (b)^{\text{\tiny $\bigstar$}} =  \bigcup_{d\mid \widetilde{m}_3} \, g(d)^{\text{\tiny $\bigstar$}}.$}
\end{itemize}
\item[(B)] If $2\mid n'$ and $2\nmid a'b'$, we have
\begin{itemize}
\item[(1)] $N_{n;a,b}(c)>0$ if and only if $g\mid c$ and $c':=\frac{c}{g}$ is an even integer satisfying $(c',a'b')=1$.
\item[(2)]  If $c\neq 0$ and $N_{n;a,b}(c)>0$, then there is a unique ideal $I\subset \mathbb Z_{n'}$ such that $c'\in I^{\text{\tiny $\bigstar$}}$. Moreover, $\mbox{lead}(I)= (c',\widetilde{m}_3)$, where $\widetilde{m}_3$ is the largest divisor of $\frac{n'}{a'b'}$ satisfying $(\widetilde{m}_3,a'b')=1\;$. 
\item[(3)]  {$\displaystyle (a)^{\text{\tiny $\bigstar$}} + (b)^{\text{\tiny $\bigstar$}} =  \bigcup_{d\mid \widetilde{m}_3,\, d \; even} \, g(d)^{\text{\tiny $\bigstar$}}.$}
\end{itemize}
\end{itemize}
\end{theorem}

%
%
%


\section{Some technical preliminaries}  \label{technical}

%
%
%

As preliminary results to the calculations in the subsequent section, we prove several identities, being somewhat charming of their own. Here $\mu$ is the Moebius function and we apply some of its most basic properties (cf.~\cite{HW}, chapters 16.3--16.4). 

\begin{lemma}  \label{klotz}
Let $n$, $d$ and $r$ be positive integers satisfying $(r,d)=1$. Then 
\[
f(n)=f_{d,r}(n) := \#\{1\leq y \leq n: \; (dy+r,n)=1\} = 
n \prod_{\genfrac{}{}{0pt}{1}{p\in \mathbb P}{p\mid n,\, p\nmid d}}
\left(1-\frac{1}{p}\right).
\] 
In particular, there is always some $y$ such that $(dy+r,n)=1$.
\end{lemma}

{\sc Proof. }
We have
\begin{align}  \label{lem0}
\begin{split}
f_{d,r}(n) &= \sum_{\genfrac{}{}{0pt}{1}{y=1}{(dy+r,n)=1}}^n 1 
= \sum_{y=1}^n \sum_{\genfrac{}{}{0pt}{1}{g\mid n}{g\mid (dy+r)}} \mu(g)
= \sum_{g\mid n} \mu(g)  \sum_{\genfrac{}{}{0pt}{1}{y=1}{dy+r\equiv 0 \bmod g}}^n 1\\
&= n \sum_{g\mid n} \frac{\mu(g)}{g}  \sum_{\genfrac{}{}{0pt}{1}{y=1}{dy+r\equiv 0 \bmod g}}^g 1, 
\end{split}
\end{align}
where the inner sum vanishes if $(d,g)\nmid r$, and equals $(d,g)$ in case $(d,g)\mid r$. Since $(r,d)=1$ by assumption, $(d,g)\mid r$ is satisfied if and only if $(d,g)=1$, and in this case the inner sum equals $1$ while it vanishes otherwise. Hence \eqref{lem0} implies
\begin{equation}\label{lem1}
f_{d,r}(n)=n \sum_{\genfrac{}{}{0pt}{1}{g\mid n}{(g,d)=1}} \frac{\mu(g)}{g}.
\end{equation}
It follows for coprime integers $m$ and $n$ that
\[
f_{d,r}(mn)= mn \sum_{\genfrac{}{}{0pt}{1}{g_1\mid m,\, g_2\mid n}{(g_1g_2,d)=1}} \frac{\mu(g_1g_2)}{g_1g_2} = f_{d,r}(m)f_{d,r}(n),
\]
i.e.~for fixed $d$ and $r$ the function $f_{d,r}$ ist multiplicative. For a prime power $p^s$, $f$ is easily evaluated by \eqref{lem1}, and we obtain
\[
f(p^s) = p^s \sum_{\genfrac{}{}{0pt}{1}{j=0}{(p^j,d)=1}}^s \frac{\mu(p^j)}{p^j}
= \left\{ \begin{array}{cl}
           p^s  &  \mbox{ if $p\mid d$, } \\
           p^s(1-\frac{1}{p}) & \mbox{ if $p\nmid d$.}
           \end{array}
           \right.
\]
Now the multiplicativity of $f$ completes the proof.

\qed

We shall now show that \eqref{atomunit} and \eqref{atom} are both representations of {\rm atom}(a), and that \eqref{atom} yields a unique representation of the elements of \eqref{atom}.
\begin{corollary}  \label{corunit}
Let $n$ be a positive integer, and let $a\in \mathbb Z_n$. Then 
${\rm atom}(a) = (a)^{\text{\tiny $\bigstar$}}$, 
and each $b\in {\rm atom}(a)$ has a unique representation $b=ax$ with $1\leq x \leq {\rm ord}(a)$, $(x, {\rm ord}(a))=1$, hence $\#\,{\rm atom}(a) = \#(a)^{\text{\tiny $\bigstar$}}=\varphi({\rm ord}(a))$. 
\end{corollary}

{\sc Proof. } 
Let $b\in a\mathbb Z_n^*$, i.e.~$b=ax$ for some $x$ satisfying $(x,n)=1$. 
Since $a(x+\ell\cdot{\rm ord}(a)) = ax +\ell(a\cdot{\rm ord}(a)) =ax$ for each integer $\ell\in \mathbb Z$, we may assume that $1\leq x \leq {\rm ord}(a)$. Since ${\rm ord}(a) \mid n$ and $(x,n)=1$, we also have $(x, {\rm ord}(a))=1$, thus $b=ax\in (a)^{\text{\tiny $\bigstar$}}$.

Conversely, let $b=ax \in (a)^{\text{\tiny $\bigstar$}}$, thus $(x, {\rm ord}(a))=1$. By Lemma \ref{klotz} there is some $y$ such that $({\rm ord}(a)\cdot y +x, n)=1$. Since $a({\rm ord}(a)\cdot y +x) = ax = b$, we have $b\in a\mathbb Z_n^*$. Up to now, we have shown that ${\rm atom}(a) =(a)^{\text{\tiny $\bigstar$}}$.

Since $ax=ax'$ with $1\leq x \leq x' \leq {\rm ord}(a)=\frac{n}{(a,n)}$ implies $x\equiv x' \bmod \frac{n}{(a,n)}$, the representation in $(a)^{\text{\tiny $\bigstar$}}$ is unique as desired. Hence  $\#\,{\rm atom}(a) = \#(a)^{\text{\tiny $\bigstar$}}=\varphi({\rm ord}(a))$. 

\qed

\begin{lemma}   \label{lemma}
Let $m$ and $k$ be positive integers. Then
\[ 
T(m,k) :=   \mathop{\sum_{d\mid m}\sum_{e\mid m}}\limits_{(d,e)=k} \frac{\mu(d)}{d}\,\frac{\mu(e)}{e}
= \frac{\vert\mu(k)\vert}{k^2} \prod_{\genfrac{}{}{0pt}{1}{p\in \mathbb P}{ p\mid m,\, p\nmid k}}
\left(1-\frac{2}{p}\right)
\]
in case $k\mid m$, and $T(m,k)=0$ if $k\nmid m$.
\end{lemma}

{\sc Proof. } 
If $k\nmid m$, then $T(m,k)=0$ is obvious. In the sequel we assume that $k\mid m$. Let $r$ and $s$ be coprime positive integers, and let 
$u$ and $v$ be defined according to $k=uv$ such that $u\mid r$ and $v\mid s$. Then
\begin{align}  \label{mult}
\begin{split}
T(rs,uv) &= 
\mathop{\sum_{d\mid r,\, d'\mid s}\sum_{e\mid r,\,e'\mid s}}\limits_{(dd',ee')=uv} \frac{\mu(dd')}{dd'}\,\frac{\mu(ee')}{ee'}  \\
&=
 \Bigg(\mathop{\sum_{d\mid r}\sum_{e\mid r}}\limits_{(d,e)=u} \frac{\mu(d)}{d}\,\frac{\mu(e)}{e}\Bigg)
  \Bigg(\mathop{\sum_{d'\mid r}\sum_{e'\mid s}}\limits_{(d',e')=v} \frac{\mu(d')}{d'}\,\frac{\mu(e')}{e'}\Bigg)
  = T(r,u)T(s,v),
\end{split}  
\end{align}
where we have used that $(dd',ee')=(d,e)(d',e')$. This multiplicativity property of $T(m,k)$ reduces the problem to the calculation of $T(m,k)$ for prime powers $m=p^n$ and $k=p^{\ell}$, say. 
Since Lemma \ref{lemma} trivially holds for $m=k=1$, we may assume $n\geq 1$. By definition we obtain
\[
T(p^n,p^{\ell}) = \mathop{\sum_{i=0}^n\sum_{j=0}^n}\limits_{\min\{i,j\}=\ell} \frac{\mu(p^i)}{p^i}\,\frac{\mu(p^j)}{p^j}.
\]
For $\ell\geq 2$ all summands apparently vanish, hence $T(p^n,p^{\ell})=0$. It remains to study the cases $\ell=0$ and $\ell =1$, where we have $T(p^n,1)= 1-\frac{2}{p}$ and $T(p^n,p)= \frac{1}{p^2}$.
By (\ref{mult}), this means that $T(m,k)=0$ if $k$ is non-squarefree. For squarefree $k$ we finally get
\[ T(m,k) = \prod_{\genfrac{}{}{0pt}{1}{p\in \mathbb P}{ p\mid k}}
\frac{1}{p^2} \prod_{\genfrac{}{}{0pt}{1}{p\in \mathbb P}{ p\mid m,\, p\nmid k}}
\left(1-\frac{2}{p}\right) =  \frac{1}{k^2} \prod_{\genfrac{}{}{0pt}{1}{p\in \mathbb P}{ p\mid m,\, p\nmid k}} \left(1-\frac{2}{p}\right).
\]
\qed
\medskip
\newline
{\bf Remark. }
Lemma \ref{lemma} is complemented by a nice identity for the Moebius function, namely
\[
\mu(k) = \mathop{\sum_{d\mid m}\sum_{e\mid m}}\limits_{[d,e]=k} \frac{\mu(d)}{d}\,\frac{\mu(e)}{e},
\]
for $k\mid m$, where $[d,e]$ denotes the least common multiple of $d$ and $e$. This result can be shown in exactly the same manner as Lemma \ref{lemma}.
\bigskip

%

Another technical tool is
\begin{lemma}   \label{lemma2}
Let $m$ and $k$ be positive integers. Then
\[ 
Q(m,k) :=   \sum_{\genfrac{}{}{0pt}{1}{d\mid m}{d \mid k}}
\frac{\vert\mu(d)\vert}{d} \prod_{\genfrac{}{}{0pt}{1}{p\in \mathbb P}{ p\mid m,\, p\nmid d}}
\left(1-\frac{2}{p}\right) =  \frac{\varphi^*(m,k)}{m}\,.
\]
\end{lemma}

{\sc Proof. }
Due to the factor $\vert\mu(d)\vert$ the term $Q(m,k)$ depends only on the squarefree kernels ${\rm rad}(m)$ of $m$ and ${\rm rad}(k)$ of $k$, and the same is apparently true for $\frac{\varphi^*(m,k)}{m}$. Therefore, we may assume without loss of generality that $m$ and $k$ are squarefree. Moreover, any prime factors of $k$ which do not divide $m$ are irrelevant. For this reason, we can also assume that $k\mid m$. Then we have
\begin{align*}
Q(m,k) &=   \sum_{d\mid k}
\frac{1}{d} \prod_{\genfrac{}{}{0pt}{1}{p\in \mathbb P}{ p\mid m,\, p\nmid d}}
\left(1-\frac{2}{p}\right)  \\
&=  \sum_{d\mid k} \Bigg(\prod_{p\in \mathbb P,\,p\mid d}\frac{1}{p}\Bigg)
\Bigg(\prod_{\genfrac{}{}{0pt}{1}{p\in \mathbb P}{ p\mid k,\, p\nmid d}}
\left(1-\frac{2}{p}\right)  
\prod_{\genfrac{}{}{0pt}{1}{p\in \mathbb P}{ p\mid m,\, p\nmid k}}
\left(1-\frac{2}{p}\right) \Bigg) \\
&=  \sum_{d\mid k} \Bigg(\prod_{p\in \mathbb P,\, p\mid d}\frac{1}{p}
\prod_{\genfrac{}{}{0pt}{1}{p\in \mathbb P}{ p\mid k,\, p\nmid d}}
\left(1-\frac{2}{p}\right) \Bigg) \Bigg( \prod_{\genfrac{}{}{0pt}{1}{p\in \mathbb P}{ p\mid m,\, p\nmid k}} \left(1-\frac{2}{p}\right)   \Bigg) \\
&= \prod_{\genfrac{}{}{0pt}{1}{p\in \mathbb P}{ p\mid k}} \left(\frac{1}{p} + \Big(1-\frac{2}{p}\Big)\right) \prod_{\genfrac{}{}{0pt}{1}{p\in \mathbb P}{ p\mid m,\, p\nmid k}} \left(1-\frac{2}{p}\right)  = \frac{\varphi^*(m,k)}{m}\; .
\end{align*}
\qed
\medskip


\section{The number of representations}

Given two ideals $I$ and $J$ in $\mathbb Z_n$ and some $c\in \mathbb Z_n$, we wish to determine the number of representations $c=u+v$  with $u\in I^{\text{\tiny $\bigstar$}}$ and $v\in J^{\text{\tiny $\bigstar$}}$. By the Reduction Lemma \ref{redlemma} we may assume that $a:={\rm lead}(I)$ and $b:={\rm lead}(J)$ are coprime.

\begin{proposition}  \label{prop1}
Let $n$ be a positive integer with coprime divisors $a$ and $b$, and let $c\in \mathbb Z_n$.
\begin{itemize}
\item[(i)] If $(c,ab)>1$, then $N_{n;a,b}(c)=0$.
\item[(ii)] Let $(c,ab)=1$. Then $m:=\frac{n}{ab}$ is a positive integer, and writing ${\rm rad}(m) = m_1 m_2m_3$ with $m_1\mid a$, $m_2\mid b$ and $(m_3,ab)=1$, we have 
\begin{align*}
N_{n;a,b}(c) &=  \frac{m}{{\rm rad}(m)}\, \varphi(m_1)\,\varphi(m_2)\,\varphi^*(m_3,c) \\
&= m  \prod_{\genfrac{}{}{0pt}{1}{p\in \mathbb P}{p\mid m, \,p\mid ab}}
\left(1-\frac{1}{p}\right)
\prod_{\genfrac{}{}{0pt}{1}{p\in \mathbb P}{p\mid n,\, p\nmid ab, \, p\mid c}} \left(1-\frac{1}{p}\right)
\prod_{\genfrac{}{}{0pt}{1}{p\in \mathbb P}{p\mid n,\,  p\nmid abc}} \left(1-\frac{2}{p}\right).
\end{align*}
\end{itemize}
\end{proposition}

{\sc Proof. }  
%

(i) Since $(c,ab)>1$, there is a prime $p\mid c$ such that $p\mid a$, say. Then $ax+by\equiv c \bmod n$ can only have a solution $x,y$ if $p\mid by$, which by $(a,b)=1$ implies $p\mid y$. But $x,y$ can contribute to $N(c)$ in (\ref{basic}) only if $(y,\frac{n}{b})=1$. It follows that $(y,a)=1$, because $a\mid \frac{n}{b}$. This contradiction shows that $N(c)=0$.

(ii) Consider the congruence $ax+by \equiv c\bmod n$ in (\ref{basic}) for some fixed $x$. Then the congruence is only solvable if $b=(b,n)\mid (ax-c)$, and in this case there exists a unique solution $y \bmod\, am$. Hence
\[
N(c) =  \mathop{\sum_{\genfrac{}{}{0pt}{1}{1\leq x \leq bm}{(x,bm)=1}} 
\sum_{\genfrac{}{}{0pt}{1}{1\leq y \leq am}{(y,am)=1}}}\limits_{ax+by\equiv c\bmod n} 1
= \sum_{\genfrac{}{}{0pt}{1}{1\leq x \leq bm}{(x,bm)=1}} 
\sum_{\genfrac{}{}{0pt}{1}{\genfrac{}{}{0pt}{1}{1\leq y \leq am}{(y,am)=1}}{by\equiv c-ax \bmod n}} 1
= \sum_{\genfrac{}{}{0pt}{1}{\genfrac{}{}{0pt}{1}{1\leq x \leq bm}{(x,bm)=1}}{b\mid (ax-c)}} 
\sum_{\genfrac{}{}{0pt}{1}{\genfrac{}{}{0pt}{1}{1\leq y \leq am}{(y,am)=1}}{y\equiv \frac{c-ax}{b} \bmod am}} 1.
\]
The required coprimality condition $(y,ma)=1$ is satisfied if and only if $(\frac{ax-c}{b},ma)=1$, and we obtain
\begin{align*}
N(c) &=   \sum_{\genfrac{}{}{0pt}{1}{\genfrac{}{}{0pt}{1}{1\leq x \leq bm}{(x,bm)=1}}
{\genfrac{}{}{0pt}{1}{b\mid (ax-c)}{(\frac{ax-c}{b},am)=1}}} 1 
=  \sum_{0\leq t < m}  \sum_{\genfrac{}{}{0pt}{1}{\genfrac{}{}{0pt}{1}{1\leq r \leq b}{(bt+r,bm)=1}}
{\genfrac{}{}{0pt}{1}{b\mid (a(bt+r)-c)}{(\frac{a(bt+r)-c}{b},am)=1}}} 1 
=  \sum_{0\leq t < m}  \sum_{\genfrac{}{}{0pt}{1}{\genfrac{}{}{0pt}{1}{1\leq r \leq b}{(bt+r,bm)=1}}
{\genfrac{}{}{0pt}{1}{ar\equiv c \bmod b}{(at+ \frac{ar-c}{b},am)=1}}} 1  \\
&=  \sum_{0\leq t < m}  \sum_{\genfrac{}{}{0pt}{1}{1\leq r \leq b}{ar\equiv c \bmod b}}  \Bigg(\sum_{d\mid (bt+r,bm)} \mu(d)\Bigg)\Bigg(\sum_{e\mid (at+ \frac{ar-c}{b},am)} \mu(e)\Bigg) \\
&=  \sum_{\genfrac{}{}{0pt}{1}{1\leq r \leq b}{ar\equiv c \bmod b}}  \sum_{d\mid bm} \mu(d) \sum_{e\mid am} \mu(e) \sum_{\genfrac{}{}{0pt}{1}{\genfrac{}{}{0pt}{1}{0\leq t < m}{bt\equiv -r \bmod d }}{at\equiv -\frac{ar-c}{b}\bmod e}} 1.
\end{align*}
Since $\mu(d)=0$ for any non-squarefree integer $d$, we conclude that
\[
N(c) =  \sum_{\genfrac{}{}{0pt}{1}{1\leq r \leq b}{ar\equiv c \bmod b}}  \sum_{d\mid b\cdot {\rm rad}(m)} \mu(d) \sum_{e\mid a\cdot{\rm rad}(m)} \mu(e) \sum_{\genfrac{}{}{0pt}{1}{\genfrac{}{}{0pt}{1}{0\leq t < m}{bt\equiv -r \bmod d }}{at\equiv -\frac{ar-c}{b}\bmod e}} 1.
\]
It follows from $(a,b)=1$ that the three divisors $m_1,m_2,m_3$ of ${\rm rad}(m)$ as defined in the statement of the proposition are pairwise coprime, and we obtain 
\begin{equation} \label{nc2}
N(c) = \sum_{\genfrac{}{}{0pt}{1}{1\leq r \leq b}{ar\equiv c \bmod b}} \sum_{d_1\mid m_1} \sum_{d_2\mid b} \sum_{d_3\mid m_3} \mu(d_1d_2d_3) \sum_{e_1\mid a} \sum_{e_2\mid m_2} \sum_{e_3\mid m_3} \mu(e_1e_2e_3)\hspace{-1cm}\sum_{\genfrac{}{}{0pt}{1}{\genfrac{}{}{0pt}{1}{0\leq t < m}{bt\equiv -r \bmod d_i \;\;(i=1,2,3) }}{at\equiv -\frac{ar-c}{b}\bmod e_i \;\; (i=1,2,3)}} \hspace{-1cm}1.
\end{equation}
The congruence system 
\begin{equation}  \label{congrsys}
\left.
\begin{array}{rcll}
bt & \equiv & -r &\bmod d_i \;\;(i=1,2,3) \\
at &\equiv & -\frac{ar-c}{b} &\bmod e_i \;\; (i=1,2,3)
\end{array}
\right\}
\end{equation}
can have a solution only if each single congruence is solvable, i.e.~$(d_1,b)\mid r$, $(d_2,b)\mid r$, $(d_3,b)\mid r$, $(e_1,a)\mid \frac{ar-c}{b}$,  $(e_2,a)\mid \frac{ar-c}{b}$ and $(e_3,a)\mid \frac{ar-c}{b}$. Due to the obvious divisibility and coprimality properties of $m_1,m_2,m_3$, four of these six conditions are satisfied per se. Only two conditions are necessary for the solvability of \eqref{congrsys}, namely $d_2=(d_2,b)\mid r$ and $e_1=(e_1,a)\mid \frac{ar-c}{b}$, which imply $b\equiv 0 \equiv r \bmod d_2$ and $a \equiv 0 \equiv \frac{ar-c}{b} \bmod e_1$. Therefore, the two congruences $\bmod\, d_2$ and $\bmod\, e_1$ in \eqref{congrsys} trivially hold for every $t$. All in all, the conditions $d_2\mid r$ and $e_1 \mid \frac{ar-c}{b}$ are necessary for the solvability of \eqref{congrsys} and yield its equivalence
with the congruence system
\begin{equation}    \label{congrsys2}
\left.
\begin{array}{rcll}
bt &\equiv &-r &\bmod d_1 \\
bt &\equiv &-r &\bmod d_3 \\
at &\equiv &-\frac{ar-c}{b}&\bmod e_2  \\
at &\equiv &-\frac{ar-c}{b}&\bmod e_3\ 
\end{array}\right\}
\end{equation}
Among the pairwise greatest common divisors of the moduli $d_1,d_3,e_2,e_3$ only $(d_3,e_3)$ may be greater than $1$.
By an extended version of the Chinese remainder theorem (for non-coprime moduli) (cf.~\cite{HW}, chapter 8), the congruence system (\ref{congrsys2}) is solvable if and only if $a r \equiv b\,\frac{ar-c}{b} \bmod (d_3,e_3)$, i.e.~$(d_3,e_3)\mid c$, and then there exists a unique solution $\bmod\;  d_1 e_2[d_3,e_3]$ with the least common multiple $[d_3,e_3]$ of the moduli $d_3$ and $e_3$. To sum up, system (\ref{congrsys}) is solvable if and only if $d_2\mid r$, $e_1\mid \frac{ar-c}{b}$, $(d_3,e_3)\mid c$, and under these conditions it has a unique solution $\bmod\; d_1 e_2[d_3,e_3]$. Using this in (\ref{nc2}) and applying the identity $d_3e_3=(d_3,e_3)[d_3,e_3]$, we obtain
\begin{align}    \label{final1}
\begin{split}
N(c) &= \sum_{\genfrac{}{}{0pt}{1}{1\leq r \leq b}{ar\equiv c \bmod b}} \sum_{d_1\mid m_1} \mu(d_1)
\sum_{\genfrac{}{}{0pt}{1}{d_2\mid b}{d_2 \mid r}} \mu(d_2)\sum_{\genfrac{}{}{0pt}{1}{e_1\mid a}{e_1\mid \frac{ar-c}{b}}} \mu(e_1)\sum_{e_2\mid m_2} \mu(e_2)  \;\; \times \\
&\quad\quad\quad\quad\quad\quad  \times  \mathop{\sum_{d_3\mid m_3} 
\sum_{e_3\mid m_3}}\limits_{(d_3,e_3)\mid c}
  \mu(d_3)\mu(e_3)\cdot \frac{m}{d_1e_2[d_3,e_3]}\\
&= m \Bigg(\sum_{d_1\mid m_1} \frac{\mu(d_1)}{d_1}\Bigg) \Bigg(\sum_{e_2\mid m_2} \frac{\mu(e_2)}{e_2}\Bigg) 
\Bigg( \mathop{\sum_{d_3\mid m_3} 
\sum_{e_3\mid m_3}}\limits_{(d_3,e_3)\mid c}   \frac{\mu(d_3)\mu(e_3)}{d_3e_3}(d_3,e_3)\Bigg) \times \\
&\quad\quad\quad\quad\quad\quad  \times 
\Bigg(\sum_{\genfrac{}{}{0pt}{1}{1\leq r \leq b}{ar\equiv c \bmod b}}   \sum_{d_2\mid (b,r)} \mu(d_2) \sum_{e_1\mid (a,\frac{ar-c}{b})} \mu(e_1)
\Bigg).
\end{split}
\end{align}
Lemma \ref{lemma} and Lemma \ref{lemma2} yield
\begin{align}  \label{final2}
\begin{split}
\mathop{\sum_{d_3\mid m_3} 
\sum_{e_3\mid m_3}}\limits_{(d_3,e_3)\mid c}   \frac{\mu(d_3)\mu(e_3)}{d_3e_3}(d_3,e_3) &=
\sum_{d\mid c} d \mathop{\sum_{d_3\mid m_3} 
\sum_{e_3\mid m_3}}\limits_{(d_3,e_3)= d}   \frac{\mu(d_3)\mu(e_3)}{d_3e_3} = \sum_{d\mid c} d\,T(m_3,d)\\
&= \sum_{\genfrac{}{}{0pt}{1}{d\mid c}{d\mid m_3}} \frac{\vert\mu(d)\vert}{d} \prod_{\genfrac{}{}{0pt}{1}{p\in \mathbb P}{ p\mid m_3,\, p\nmid d}}
\left(1-\frac{2}{p}\right) =  Q(m_3,c) \\
&= \frac{\varphi^*(m_3,c)}{m_3}.
\end{split}
\end{align}
Assuming $r$ to be a solution of $ar\equiv c \bmod b$, it follows from $(b,c)=1$ that $(b,r)=1$. In addition, $1=(a,c)=(a,ar-c)$ implies $(a,\frac{ar-c}{b})=1$. Hence
\begin{equation}   \label{VorschlagTorsten}
\sum_{\genfrac{}{}{0pt}{1}{1\leq r \leq b}{ar\equiv c \bmod b}}   \sum_{d_2\mid (b,r)} \mu(d_2) \sum_{e_1\mid (a,\frac{ar-c}{b})} \mu(e_1) =
\sum_{\genfrac{}{}{0pt}{1}{1\leq r \leq b}{ar\equiv c \bmod b}} \mu(1)\mu(1)
= \sum_{\genfrac{}{}{0pt}{1}{1\leq r \leq b}{ar\equiv c \bmod b}} 1 =1,
\end{equation}
since the congruence has exactly one solution $r \bmod b$.
Inserting \eqref{VorschlagTorsten} and (\ref{final2}) into (\ref{final1}), and using a standard identity for the totient function, we obtain
\begin{align*}
N(c) 
&= m \; \frac{\varphi(m_1)}{m_1} \; \frac{\varphi(m_2)}{m_2} \; \frac{\varphi^*(m_3,c)}{m_3} 
= \frac{m}{{\rm rad}(m)}\,\varphi(m_1)\varphi(m_2)\varphi^*(m_3,c)\\
&= m \prod_{\genfrac{}{}{0pt}{1}{p\in \mathbb P}{p\mid m_1}}
\left(1-\frac{1}{p}\right)
\prod_{\genfrac{}{}{0pt}{1}{p\in \mathbb P}{p\mid m_2}}
\left(1-\frac{1}{p}\right)
\prod_{\genfrac{}{}{0pt}{1}{p\in \mathbb P}{p\mid m_3, \, p\mid c}} \left(1-\frac{1}{p}\right)
\prod_{\genfrac{}{}{0pt}{1}{p\in \mathbb P}{p\mid m_3,\,  p\nmid c}} \left(1-\frac{2}{p}\right).
\end{align*}
The proof of (ii) is completed by the fact that 
\[
p\mid m_1 \Leftrightarrow (p\mid m \mbox{ and } p\mid a), 
\quad p\mid m_2 \Leftrightarrow (p\mid m \mbox{ and } p\mid b), \quad 
 p\mid m_3 \Leftrightarrow (p\mid n \mbox{ and } p\nmid ab).
\]
\qed

{\sc Proof of Theorem \ref{thm1}. }

(i)  Clearly, $g\mid a$, $g\mid b$ and $g\mid n$. If $g\nmid c$, then $ax+by\equiv c \bmod n$ has no solution and thus $N(c)=0$ by (\ref{basic}).

(ii) Since $g\mid c$, all numbers $n':=\frac{n}{g}$, $a':=\frac{a}{g}$, $b':=\frac{b}{g}$ and $c':=\frac{c}{g}$ are integers satisfying $(a',b')=1$. Hence $N_{n;a,b}(c) =  N_{n';a',b'}(c')$ by Reduction Lemma \ref{redlemma} (ii). 
As a consequence of Proposition \ref{prop1}~(i) it follows that $N_{n;a,b}(c) = 0$ if $(c',a'b')>1$.
In case $(c',a'b')=1$, the identities of \eqref{genform} follow right away from Proposition \ref{prop1}~(ii).

\qed

{\sc Proof of Corollary \ref{cor1}. }
By definition of $N_{n;a,b}(c)$, we have  $c\in (a)^{\text{\tiny $\bigstar$}} + (b)^{\text{\tiny $\bigstar$}}$ if and only if $N_{n;a,b}(c) >0$. By Theorem \ref{thm1}, the non-vanishing of $N_{n;a,b}(c)$ necessarily requires (i) and (ii). Under these two conditions, we have $N_{n;a,b}(c) >0$ unless the factor $\varphi^*(m_3,c')$ in (\ref{genform}) vanishes, which means that $2\mid m_3$ and $2\nmid c'$. Hence $\varphi^*(m_3,c')>0$ if and only if $2\nmid m_3$ or $2\mid c'$, which in turn is equivalent with (iii).

\qed

{\sc Proof of Corollary \ref{cor2}. }
Let $c:={\rm lead}(I)$. Then it suffices to show that 
\begin{equation}   \label{compare}
N_{n;a,b}(c) = N_{n;a,b}(cx) \quad\quad\quad \big(1\leq x \leq \tfrac{n}{c},\;\;  
(x, \tfrac{n}{c})=1\big). 
\end{equation}
We set $g:=(a,b)$.

\underline{Case 1:} $\;g\nmid c$.\newline
We know from Theorem \ref{thm1} (i) that $N_{n;a,b}(c)=0$. Moreover, there is a prime $p$ such that $p\mid g$, but $p\nmid c$. Since $p\mid n$, it follows that $p\mid \frac{n}{c}$, which by $(x, \frac{n}{c})=1$ implies $p\nmid x$. We obtain $p\nmid cx$, thus $g\nmid cx$. Now Theorem \ref{thm1} (i) tells us that $N_{n;a,b}(cx)=0=N_{n;a,b}(c)$.

\underline{Case 2:} $\;g\mid c$.\newline
As before, we set $n':=\frac{n}{g}$, $a':=\frac{a}{g}$, $b':=\frac{b}{g}$ and $c':=\frac{c}{g}$.
If $(c',a'b')>1$, then $(c'x,a'b')>1$, and Theorem \ref{thm1} (ii) yields $N_{n;a,b}(cx)=0=N_{n;a,b}(c)$. Hence we may assume $(c',a'b')=1$. The fact that $(x,\frac{n'}{c'})=(x,\frac{n}{c})=1$ implies $(c'x,a'b')=(x,a'b')=1$ enables us to compare $N_{n;a,b}(c)$ and $N_{n;a,b}(cx)$ by \eqref{genform}.
Since $m,m_1,m_2,m_3$ as defined in Theorem \ref{thm1} (ii) do not depend on $c$ or $cx$, respectively, \eqref{compare} would follow from the second identity of \eqref{genform} if we can prove that $\varphi^*(m_3,c')=\varphi^*(m_3,c'x)$ for all $x$ satisfying $(x,\frac{n'}{c'})=(x,\frac{n}{c})=1$.
Therefore, it suffices to show 
that $p\mid c' \Leftrightarrow p \mid c'x$ for all primes $p\mid m_3$. The direction from left to right is trivial. Conversely, we assume that $p\mid x$ and have to deduce that $p\mid c'$. It follows from $p\mid m_3$ that $p\mid n'$. Since $(x, \frac{n'}{c'})=1$ and $p\mid x$, we obtain $p\nmid \frac{n'}{c'}$. Together this indeed implies that $p\mid c'$. 

\qed

%

{\sc Proof of Theorem \ref{thm2}. }
We prove (A) and leave the similar proof of (B) to the reader, the main difference being the obvious fact that in case (B) the sumset $(a)^{\text{\tiny $\bigstar$}} + (b)^{\text{\tiny $\bigstar$}}$ contains only even integers.

(A1) This is an immediate consequence of Corollary \ref{cor1} by virtue of our condition $2\nmid n'$ or $2\mid a'b'$.

(A2) By (A1) we have $g\mid c$ and $(c',a'b')=1$. We set $d:=(c',\widetilde{m}_3)$ and $k:=\frac{c'}{d}$, hence $(k,\frac{\widetilde{m}_3}{d})=1$. Observe that $k\neq 0$. We factorise $n'=a'b'\widetilde{m}_1\widetilde{m}_2\widetilde{m}_3$ in such a way that 
$\widetilde{m}_1$ contains only prime factors of $a'$ and $\widetilde{m}_2$ contains only prime factors of $b'$. 
Take notice of the fact that $(c',a'b')=1$ implies $(c',a'b'\widetilde{m}_1\widetilde{m}_2)=1$. Thus $(k,a'b'\widetilde{m}_1\widetilde{m}_2)=1$,
which yields
\[
(k,\mbox{ord}(d))=(k,\tfrac{n'}{d})=(k,a'b'\widetilde{m}_1\widetilde{m}_2\tfrac{\widetilde{m}_3}{d}) = (k,\tfrac{\widetilde{m}_3}{d}) = 1.
\]
Therefore, $c'=dk \in (d)^{\text{\tiny $\bigstar$}}$, where $d= \mbox{lead}(I)$ for $I:=(d)$. Since the uniqueness of $I$ and its leader $d$ are clear (cf.~introduction on atoms), the proof of (A2) is complete.

(A3) It suffices to prove 
\begin{equation}   \label{help}
(a')^{\text{\tiny $\bigstar$}} + (b')^{\text{\tiny $\bigstar$}} =  \bigcup_{d\mid \widetilde{m}_3} \, (d)^{\text{\tiny $\bigstar$}},
\end{equation}
since (A3) then follows by Reduction Lemma \ref{redlemma} (i).
First assume that not both integers $a'$ and $b'$ are equal to $1$. This means that $\widetilde{m}_3\neq n'$, hence $n'$ is not a divisor of $\widetilde{m}_3$.
Since condition (ii) of Corollary \ref{cor1} is violated for $c=0$, we know that $N_{n';a',b'}(0)=0$. It follows from (A2) that 
\begin{equation}  \label{union}
0 \notin (a')^{\text{\tiny $\bigstar$}} + (b')^{\text{\tiny $\bigstar$}} \subseteq \bigcup_{d\mid \widetilde{m}_3} \, (d)^{\text{\tiny $\bigstar$}}.
\end{equation}

On the other hand, let $d$ be any divisor of $\widetilde{m}_3$ and $(k,\tfrac{n'}{d})=1$, i.e.~$dk$ is an arbitrary element of the union on the righthand side of (\ref{union}). In order to complete the proof of \eqref{help} it suffices to show that $N_{n';a',b'}(dk)>0$, which by (A1) requires $(dk,a'b')=1$. Since $d\mid \widetilde{m}_3$, we have $(d,a'b')=1$ by definition, and thus $a'b'\mid \frac{n'}{d}$. Since $(k,\frac{n'}{d})=1$, we also have $(k,a'b')=1$.

We are left with the special case $a'=b'=1$, when $N_{n';1,1}(c')>0$ for all $c'$ by Corollary \ref{cor1}. Clearly, $\widetilde{m}_3=n'$, hence by the above argument and $\{0\}=\{0\}^{\text{\tiny $\bigstar$}}$,
\begin{equation*} 
0 \in \mathbb Z_{n'}^* + \mathbb Z_{n'}^* \subseteq \{0\} \cup \bigcup_{d\mid n',\, d\neq n'} \, (d)^{\text{\tiny $\bigstar$}}
= \bigcup_{d\mid n'} \, (d)^{\text{\tiny $\bigstar$}}.
\end{equation*}
The converse inclusion follows trivially from the fact that $N_{n';1,1}(c')>0$ for all $c'$, i.e.~$\mathbb Z_{n'}^* + \mathbb Z_{n'}^* = \mathbb Z_{n'}$. Therefore, \eqref{help} holds in all cases.

\qed

%
%
%
%
%
%
%


\section{Application to Cayley graphs}\label{sec:gr}

Besides the fact that adding multiplicatively defined objects is a particularly interesting
study subject for number theorists, we now present an application in graph theory.
Cayley graphs model certain algebraic properties of groups in terms of adjacency of vertices
in graphs. Given a finite additive group $G$ and a subset $S\centernot=\emptyset$ with $-S=S$,
we define the Cayley graph $\Cay(G, S)$ as follows. The vertices are identified with the elements of $G$. Two vertices $x,y\in G$ are adjacent if and only if $x-y\in S$.
The set $S$ is called the symbol of the Cayley graph $\Cay(G, S)$. In order to avoid loops,
one usually requires $0\centernot\in S$. If $G$ is a cyclic group $\mathbb Z_n$, then we obtain
the important subclass of circulant graphs (their adjacency matrices being circulant matrices).

Among these graphs there are those with integer eigenvalues (of their respective adjacency matrices),
called the integral circulant graphs. According to \cite{WSO}, these can be characterized
as follows. Let $\mathbb Z_n = \{0,\ldots,n-1\}$ be the set of vertices and choose a subset $D$ of the positive divisors of $n$.
With each divisor $d$ of $n$ we associate a set $S_n(d)=\{x\in\mathbb Z_n:~(x,n)=d\}$. Setting
$S(D)=\bigcup_{d\in D} S_n(d)$, we obtain the integral circulant graph
$\ICG(n,D):=\Cay(\mathbb Z_n, S(D))$ with $n$ vertices and divisor set $D$.

Obviously, the sets $S_n(d)$ are nothing but the atoms of $\mathbb Z_n$. Hence
\[
   \ICG(n,D) = \Cay(\mathbb Z_n, S(D))=\Cay(\mathbb Z_n, \bigcup_{d\in D} \atom(d)) = \Cay(\mathbb Z_n, \bigcup_{d\in D} (d)^{\text{\tiny $\bigstar$}}).
\]
Since this is an actual characterization of integral circulant graphs, we see that
integrality of Cayley graphs over $\mathbb Z_n$ can be determined by whether their symbol sets
can be partitioned into complete sets of atoms. It is worth noting that this view
even extends to integral Cayley graphs over finite abelian groups in general \cite{ALP}.

Due to the property of $\ICG(n,D)$ being a circulant graph, the neighbourhood of every vertex 
looks basically the same (except for a translation), just note that $x-y\in S$ if and
only if $(x+s)-(y+s)\in S$. So if we want to explore the neighbourhoods of the vertices
of some graph $\ICG(n,D)$, we can restrict ourselves to vertex $0$. Clearly,
the neighbours of vertex $0$ are given by the set $S(D)$. In order to explore the neighbourhood
of the neighbourhood of vertex $0$, we need to form the set $S(D)+S(D)$. 
This is because in a Cayley graph the act of moving from some vertex to one of its neighbours is
the same as adding some element from its symbol set.

Since $S(D)$ is a disjoint union of atoms, we just need to determine all sumsets of pairs of those atoms. By Theorem \ref{thm2} we know that these sumsets are again disjoint unions of atoms.
This gives us some interesting information on how vertices are visited when exploring an
integral circulant graph. Not only does the neighbourhood of some neighbourhood decompose into
complete atoms, but we can also tell how often the vertices of each atom get discovered.

Forming the overall union of the atoms that make up the neighbourhood of the neighbourhood of 
vertex $0$, we effectively determine all vertices in $\ICG(n,D)$ whose distance from vertex $0$
is at most two. Filtering out the vertex $0$ at distance level $0$ (corresponding to the
atom $(0)^{\text{\tiny $\bigstar$}}$ and the neighbourhood of vertex $0$ at distance level $1$, the remaining 
atoms form the distance level $2$. We can continue this exploration process until we have
discovered every single vertex of the graph (in which case the distance level equals the diameter
of the graph).

Recording the adjacencies of vertex $0$ with the vertices of one or more distance levels 
and extending this into a circulant adjacency matrix, we obtain a generalized distance
matrix of the graph $\ICG(n,D)$. By construction, each such matrix represents an
integral circulant graph for some particular divisor set and the divisor sets of all
distance level graphs form a partition of the divisor set of $n$. 
Choosing the consecutive distances $1,\ldots,r$ for some $1\leq r\leq n$, we obtain 
the so-called distance powers of $\ICG(n,D)$.

As an example, consider the graph
\[
   \Gamma:=\ICG(60,\{3,10\})=\Cay(\mathbb Z_{60}, \atom(3)\cup\atom(10)).
\]
We shall determine the divisor set of its second distance power
\[
   \Gamma^{(2)}=\Cay\left(\mathbb Z_{60}, \bigcup\limits_{u,v\in\{3,10\}}\left( \atom(u)+\atom(v)\right)\setminus \{0\}\right).
\]
Let us evaluate the involved atom sums:
\[
\begin{split}
(3)^{\text{\tiny $\bigstar$}}+(3)^{\text{\tiny $\bigstar$}}   &= \bigcup\limits_{2\mid d,\, d\mid 20} 3(d)^{\text{\tiny $\bigstar$}} = S_{60}(\{0,6,12,30\}),\\
(3)^{\text{\tiny $\bigstar$}}+(10)^{\text{\tiny $\bigstar$}}  &= \bigcup\limits_{d\mid 1} ~(d)^{\text{\tiny $\bigstar$}} = S_{60}(\{1\}),\\
(10)^{\text{\tiny $\bigstar$}}+(10)^{\text{\tiny $\bigstar$}} &= \bigcup\limits_{2\mid d,\, d\mid 6} 10(d)^{\text{\tiny $\bigstar$}} = S_{60}(\{0,20\}).
\end{split}
\]
Hence it follows that
\[
  \Gamma^{(2)}=\ICG(60,\{1,6,12,20,30\}).
\]
Noting the atom leaders of the respective distance levels, we additionally obtain the information that during the exploration process of the distance levels,
none of the atom vertex sets has been explored more than once, with exception of $\{0\}=\atom(0)$:
$$
\begin{array}{ll}
\text{Level 0:} & 0 \\
\text{Level 1:} & \underbrace{3,9,21,27,33,39,51,57}_{\atom(3)},\underbrace{10,50}_{\atom(10)} \\
\text{Level 2:} & 0,\underbrace{1,7,11,13,17,19,23,29,31,37,41,43,47,49,53,59}_{\atom(1)},\underbrace{6,18,42,54}_{\atom(6)},\\
 & \underbrace{12,24,36,48}_{\atom(12)},\underbrace{20,40}_{\atom(20)},\underbrace{30\phantom{,}}_{\atom(30)}
\end{array}
$$

\end{document}